\documentclass[10pt,reqno]{amsart}
\usepackage{amsmath,amssymb}

\usepackage{latexsym}

\newtheorem{thm}{Theorem}
\newtheorem{prop}{Proposition}
\newtheorem{lem}[prop]{Lemma}
\newtheorem{nota}[prop]{Notation}

\newtheorem{defi}[prop]{Definition}
\newtheorem{rem}[prop]{Remark}

\newenvironment{pf}{\begin{trivlist}\item[]{\bf Proof.\ }}
{\mbox{}\hfill\rule{.08in}{.08in}\end{trivlist}}

\def\R{\mathbb{R}}
\def\N{\mathbb{N}}

\title{Bounded geometry, growth and topology}

\thanks{Work supported by the Project
``Internazionalizzazione" {\em {Propriet\`a asintotiche di
variet\`a e di gruppi discreti}} of MIUR of Italy.}

\author{Renata Grimaldi}

\address{Renata Grimaldi: Universit\`a Degli Studi di Palermo, Dipartimento di Metodi e Modelli Matematici,
Viale delle Scienze 90128 Palermo, Italy}
\email{grimaldi@unipa.it}

\author{Pierre Pansu}
\address{Pierre Pansu: \'Ecole Normale Sup\'erieure, DMA-ENS
45 rue d'Ulm,  F-75230 Paris Cedex 05, France}
\email{Pierre.Pansu@ens.fr}

\begin{document}

\begin{abstract}
We characterize functions which are growth types of Riemannian manifolds of bounded geometry.

\vspace{0.1cm} \noindent {\bf Keywords:} Bounded geometry, growth
types, finite topological type, graphs, quasi-isometries.

\vspace{0.1cm} \noindent {\bf MSC Subject:} 53C20. 
\end{abstract}

\maketitle

\section{Introduction and results}

In this paper we will be mostly interested in manifolds of
bounded geometry. Such spaces arise naturally when one deals with
non-compact Riemannian manifolds, for example universal coverings
of compact manifolds lie within this class of open manifolds.
Roughly speaking, a manifold of bounded geometry can be seen as a
non-compact manifold whose geometric complexity is bounded.

Our aim is to understand what are the possible growth types of
connected Riemannian manifolds of bounded geometry, continuing work by M.
Badura, \cite{Ba}.

Recall that two nondecreasing functions $v$, $w:\N\to\R_{+}$ have the same growth type if there exists an integer $A$ such that for all $n\in \N$,
\begin{eqnarray*}
w(n)\leq Av(An+A)+A,\quad v(n)\leq Aw(An+A)+A.
\end{eqnarray*}
The {\em growth type} of a connected Riemannian manifold $M$ is the growth type of the volume of balls, $n\mapsto vol(B(o,n))$. This does not depend on the choice of origin $o$.

We shall use some notions and results of the papers \cite{FG1, FG2, G1} by
L. Funar and R. Grimaldi.

\begin{defi}
A non-compact Riemannian manifold $(M,g)$ has {\em bounded
geometry} if the sectional curvature $K$ and the injectivity
radius $i_g$ satisfy $$|K|\leq 1, \ \  i_g \geq 1.$$
\end{defi}

\begin{defi}
A non-compact manifold $M$ is of {\em finite topological type}
if $M$ admits an exhaustion by compact submanifolds $M_i$ with $\partial
M_i$ all diffeomorphic to a fixed manifold $V_0$.
\end{defi}

\begin{defi}
A function $v:\N\to\R_+$ has \emph{bounded growth of derivative}
(abbreviation: $v$ is a \emph{bgd-function}) if there exists a
positive constant $L$ such that, for all $n\in\N$,
\begin{eqnarray*}
\frac{1}{L}\leq v(n+2)-v(n+1)\leq L(v(n+1)-v(n)).
\end{eqnarray*}
\end{defi}

The following statement follows from Bishop-Gromov's inequality.

\begin{prop}
\label{necessary}
Let $M$ be a connected Riemannian manifold of bounded geometry. Then the
growth function $n\mapsto vol(B(o,n))$ is a bgd-function.
\end{prop}

The main result of this paper is a converse to this statement.

\begin{thm}
\label{infinite}
Let $M$ be a connected manifold.
\begin{enumerate}
  \item If $M$ has finite topological type, every bgd-function belongs to the growth type of a Riemannian manifold of bounded geometry diffeomorphic to $M$.
  \item If $M$ has infinite topological type, a bgd-function $v$ belongs to the growth type of a Riemannian manifold of bounded geometry diffeomorphic to $M$ if and only if
$$\lim_{n\to\infty} \frac{v(n)}{n}=+\infty.$$
\end{enumerate}
\end{thm}

To get a complete characterization of growth types of Riemannian manifolds of bounded geometry, one would need a neat criterion for a growth type to contain a bgd-function. We leave this as an open question.

The proof of Theorem \ref{infinite} consists in constructing trees with prescribed growth, and plumbing Riemannian manifolds with boundary according to the combinatorial scheme provided by these trees. The pieces are provided by an exhaustion of the given manifold. Therefore their geometries are essentially unknown. Nevertheless, one can arrange so that these geometries do not interfere much with growth.

A more detailed sketch of the proof in given in subsection \ref{sketch}.

\section{Necessary conditions}

\subsection{Bounded growth of derivative}

Here we prove Proposition \ref{necessary}.

Since $M$ is connected, there exists a point $x$ at distance
$n+\frac{3}{2}$ from $o$. Then $B(o,n+2)\setminus B(o,n+1)$ contains
$B(x,\frac{1}{2})$ whose volume is bounded below in terms of sectional
curvature and injectivity radius.
This provides us with the required lower bound on $vol(B(o,n+2))-vol(B(o,n+1))$.

The upper bound follows from Bishop-Gromov's inequality in the following form. Let $M$ be a
complete $m$-dimensional Riemannian manifold with Ricci curvature $\geq
-(m-1)\kappa^2$. Let
$o\in M$. Then
\begin{eqnarray*}
r\mapsto \frac{vol(\partial B(o,r))}{\sinh(\kappa r)^{m-1}}
\end{eqnarray*}
is nonincreasing.

This implies that for all $r>0$,
\begin{eqnarray*}
vol(\partial B(o,r+1))&=&\sinh(\kappa (r+1))^{m-1}\frac{vol(\partial
  B(o,r+1))}{\sinh(\kappa (r+1))^{m-1}}\\
&\leq&\sinh(\kappa (r+1))^{m-1}\frac{vol(\partial B(o,r))}{\sinh(\kappa r)^{m-1}}.
\end{eqnarray*}
If $r\geq 1$, $\sinh(\kappa r)\geq e^{\kappa r}(1-e^{-2\kappa})/2$ and
$\sinh(\kappa (r+1))\leq e^{\kappa(r+1)}/2$, thus
\begin{eqnarray*}
\frac{\sinh(\kappa (r+1))}{\sinh(\kappa r)}\leq\frac{e^{\kappa}}{1-e^{-2\kappa}},
\end{eqnarray*}
leading to
\begin{eqnarray*}
vol(\partial B(o,r+1))\leq L\,vol(\partial B(o,r)),
\end{eqnarray*}
with $L=(e^{\kappa}/1-e^{-2\kappa})^{m-1}$.
Integrating from $n$ to $n+1$ yields
\begin{eqnarray*}
vol(B(o,n+2))-vol(B(o,n+1))\leq L(vol(B(o,n+1))-vol(B(o,n)),
\end{eqnarray*}
provided $n\geq 1$.

\subsection{A variant}

For future use, let us state the following easy variant of Bishop-Gromov's inequality.
\begin{lem}
\label{BG}
Let $M$ be a complete $m$-dimensional Riemannian manifold with totally geodesic boundary and Ricci curvature $\geq -(m-1)\kappa^2$. Let $C$ be an open and closed subset of the boundary. Let $U_r$ denotes its tubular neighborhood of width $r$. Then
\begin{eqnarray*}
r\mapsto \frac{vol(\partial U_r)}{\cosh(\kappa r)^{m-1}}
\end{eqnarray*}
is nonincreasing.
\end{lem}

As above, this implies that if $\kappa$ is small enough, then for all $k\geq 0$,
\begin{eqnarray*}
vol(U_{k+2} \setminus U_{k+1})\leq 2\,vol(U_{k+1} \setminus U_{k}).
\end{eqnarray*}

\subsection{Finite topological type}

\begin{prop}
\label{necessaryfinite}
Let $M$ be a connected Riemannian manifold of bounded geometry. Let $v(n)=vol(B(o,n))$ denote its volume growth.
Assume that $v(n)/n$ does not tend to $+\infty$. Then $M$ has finite topological type.
\end{prop}

\begin{pf}
This follows from the proof of the Funar-Grimaldi theorem, \cite{FG2}. In that paper, the first step in the proof shows that given a
function $v$ with linear growth, there exists a constant $c$ and a
sequence $n_j$ tending to $+\infty$ such that for all $j$, $v(n_j
+1)-v(n_j)\leq c$. The rest of the proof does not use linear
growth any more. Thus the proof works under the weaker assumption
that $v(n+1)-v(n)$ does not tend to $+\infty$. This assumption holds if
$v(n)/n$ does not tend to $+\infty$.
\end{pf}

\section{Sufficient conditions}

\subsection{Scheme of the construction}
\label{sketch}

A manifold diffeomorphic to $M$ will be obtained by gluing together pieces
according to the pattern given by an admissible rooted tree $T$.

\begin{defi}
\label{deftree}
Say a rooted tree $T$ is {\em admissible} if the following holds.
\begin{itemize}
  \item Each vertex of $T$ has either zero, one or two children.
  \item There is a ray (subtree homeomorphic to a half-line) emanating from the root, called the {\em trunk}, which plays a special role: the component of the root in the complement in the tree of any edge of the trunk is a finite tree.
\end{itemize}
\end{defi}

All pieces are compact Riemannian manifolds with boundary with bounded
geometry. The metric is a product in a neighborhood of the boundary. When disconnected, the boundary is split into two open and closed parts $\partial^{-}$ and $\partial^{-}$. Here is the catalog where pieces will be picked.
\begin{enumerate}
  \item A given sequence $Q_j$ of compact Riemannian manifolds with boundary with bounded geometry. $\partial^{+}Q_{j}$ is assumed to be isometric to $\partial^{-}Q_{j+1}$.
  \item For each $j$, a Riemannian manifold $R_j$ diffeomorphic to a product $\partial^{+}Q_{j}\times[0,1]$ with a disk removed, with $\partial_{-}R_{j}$ isometric to $\partial^{+}Q_{j}$ and $\partial_{+}R_{j}$ isometric to a disjoint union $S^{m-1}\cup\partial^{-}R_{j}$.
  \item Cylinder $K=S^{m-1}\times[0,\ell]$.
  \item Half sphere $HS=S^{m}_{+}$.
  \item Join $J$, diffeomorphic to a sphere with $3$ balls removed, with $\partial^{-}J=S^{m-1}$ a round sphere, and $\partial^{+}J=S^{m-1}\cup S^{m-1}$ a disjoint union of two round spheres.
  \end{enumerate}
Here are rules for the lego game. Let $S$ denote the set of vertices of the trunk having exactly one child. Then $S$ is a union of intervals $x_{n_j},x_{n_j +1},\ldots,x_{n_j +t_j -1}$ of lengths $t_j$.
\begin{enumerate}
  \item A half-sphere is chosen for the root vertex.
  \item A piece $Q_j$ is affected collectively to the vertices of the segment $[x_{n_j},\ldots,$ $x_{n_j +t_j -1}]$ of the trunk. For $n_j +t_j \leq k<n_{j+1}$, the vertex $x_{k}$ of the trunk is equipped with $R_j$.
  \item For non trunk vertices, joins, cylinders or half-spheres are chosen depending wether the number of children is $2$, $1$ or $0$.
\end{enumerate}

Lemma \ref{diffeo} asserts that the diffeomorphism type of the resulting manifold $R_{T}$ does not depend on the choice of tree. Lemma \ref{tree} shows how to construct an admissible tree $T_v$ adapted to a
prescribed function $v$. The required pieces are constructed in Proposition \ref{pieces}. Then Proposition \ref{prall} asserts that an integer valued simplification of the growth function of $R_{T_{v}}$ is equivalent to $v$. Meanwhile, one encounters twice the need to change representative of a growth type to improve its properties, Lemmas \ref{w} and \ref{suplinear}. The proof of Theorem \ref{infinite} is completed in subsection \ref{proofinfinite}.

\subsection{Matching diffeomorphism types}

\begin{lem}
\label{diffeo}
Let $T$ be an admissible rooted tree. Glue pieces according to the above three rules. The the diffeomorphism type of resulting manifold $R_T$ does not depend on $T$, only on the sequence $Q_j$.
\end{lem}

\begin{pf}
Let $T$ be an admissible rooted tree. Let
\begin{eqnarray*}
S=\bigcup_{j}[n_j ,n_j +t_j -1]
\end{eqnarray*}
be the set of (indices of) single child trunk vertices in $T$. Let $T'$ be the tree obtained from a ray $\{x_0 ,x_1 ,\ldots\}$ by adding an edge emanating from $x_n$ if and only if $n\notin S$. This is again an admissible tree. We claim that $R_{T}$ and $R_{T'}$ are diffeomorphic.

Cut $T$ (resp. $T'$) at the edge $[x_{n_j},x_{n_{j}+1}]$. By definition of admissibility, this results in finite trees, and there are corresponding manifolds with boundary $S_j$ and $S'_j$, whose boundaries are diffeomorphic to $\partial^{+}Q_{j}$. Then $S_j$ is diffeomorphic to the connected sum of $S'_{j}$ with a finite number of spheres, i.e. to $S'_j$. As $j$ increases, one can arrange that the diffeomorphism $S_{j+1}\to S'_{j+1}$ extends the previous diffeomorphism $S_{j}\to S'_{j}$, and in the limit, one gets a diffeomorphism $R_{T}\to R_{T'}$.
\end{pf}

\begin{rem}
\label{diff}
Every connected non compact manifold is diffeomorphic to some $R_T$.
\end{rem}
Indeed, let $M_j$ be an exhaustion of $M$ by connected compact submanifolds with boundary, such that $M_0$ is a disk. As we shall see in subsection \ref{req}, one can easily construct a bounded geometry metric on $M$ which is a product in a neighborhood of each $\partial M_{j}$. Letting $Q_{j}=M_{j+1}\setminus M_{j}$, one can construct $R_{j}$ as well. Inserting $R_{j}$ capped with a half sphere between $Q_{j}$ and $Q_{j+1}$ does not change the diffeomorphism type. The resulting manifold is $R_{T}$ where $T$ is the admissible tree obtained by adding an edge to every second vertex of a ray.

\subsection{Admissible trees}
\label{trees}

\begin{lem}
\label{tree}
Let $v:\N\to\N$ satisfy
\begin{itemize}
\item $v(0)=1$.
\item for all $n\in\N$, $2\leq v(n+2)-v(n+1)\leq 2(v(n+1)-v(n))$.
\item $v(n)=O(\lambda^n)$ for some $\lambda<2$.
\end{itemize}
Fix a subset $S\subset \N$ of vanishing lower density, i.e.
\begin{eqnarray*}
\liminf_{n\to \infty}\frac{|S\cap\{0,\ldots,n\}|}{n}=0.
\end{eqnarray*}
There exists an admissible rooted tree $T_{S,v}$ with bounded geometry and with growth exactly $v$ at the root.
\end{lem}

\begin{pf}
At the same time as we construct the tree inductively, we choose an ordering on the children of each vertex, and order vertices lexicographically. Put $v(1)-v(0)$ edges at the root. Assume the tree has been constructed up to level $n$. Since $v(n+1)-v(n)\leq 2(v(n)-v(n-1))$, one can glue a
total of $v(n+1)-v(n)$ edges to the $v(n)-v(n-1)$ vertices at distance
$n$ in such a way that
\begin{itemize}
  \item each vertex receives at most $2$ edges,
  \item the first one in lexicographical order receives $1$ or $2$ edges depending wether $n\in S$ or not,
  \item a maximum of them get none, and preferably the last ones in lexicographical order.
\end{itemize}
Since for all $n$, $v(n+1)-v(n)\geq 1$, the resulting graph
is connected. In fact, it is a tree with valency $\leq 3$. The trunk consists of one vertex at each level, the smallest in lexicographical order. Let us denote them by $x_k$, $k\in\N$.

Let $e=[x_{k},x_{k+1}]$ be an edge of the trunk. Assume that there exists an infinite ray emanating from the root and avoiding $e$. Let $o=y_0 ,y_1 ,\ldots$ denote its consecutive vertices. Then $y_{k+1} \not=x_{k+1}$. Since $y_{k+1}$ has at least one child, our construction forces $x_{k+1}$ to have exactly $2$ children, all of which come before $y_{k+2}$ in lexicographic order, unless $k+1\in S$. Since $y_{k+2}$ has at least one child, both of $x_{k+1}$'s children have exactly $2$ children, unless $k+2\in S$. And so on. Consider the tree obtained from the subtree emanating from $x_{k+1}$ by collapsing all edges $[x_{n},x_{n+1}]$ for $n\in S$, $n\geq k$. This is a regular binary rooted tree. This gives a lower bound of $2^{n-s(n)-k-1}$ for the number of vertices at level $n$ in $T$, where $s(n)$ denotes the number of elements of $S$ in $\{0,\ldots,n\}$. Since, by assumption, $s(n)/n$ takes arbitrarily small values, this contradicts the fact that $v(n)=O(\lambda^n)$ for some $\lambda<2$.
\end{pf}

The assumptions in Lemma \ref{tree} are not restrictive, as the following Lemma shows.
\begin{lem}
\label{w}
Let $v$ a bgd-function. Then there exists an integer valued nondecreasing function $w$ having the same growth type as $v$ such that
\begin{itemize}
\item $w(0)=1$.
\item for all $n\in\N$, $2\leq w(n+2)-w(n+1)\leq 2(w(n+1)-w(n))$.
\item $w(n)=O(\lambda^n)$ for some $\lambda<2$.
\end{itemize}
If furthermore $v(n+1)-v(n)$ tends to $+\infty$, so does $w(n+1)-w(n)$.
\end{lem}

\begin{pf}
Let $L$ be the constant controlling the growth of the derivative of $v$. Let $\ell$ be an integer such that $\ell>\log_{2}(L)$. Define a new function $z$ at multiples of $\ell$ by $z(k\ell)=v(k)$, and extend $z$ recursively at other integers as follows.
\begin{eqnarray*}
z(k\ell+s+1)=z(k\ell+s)+L^{s/\ell}\frac{L^{1/\ell}-1}{L-1}(v(k+1)-v(k)).
\end{eqnarray*}
This formula, which, when summing up, implies that $z((k+1)\ell)=z(k\ell)+v(k+1)-v(k)$, is compatible with the previous definition. For $k\ell\leq n\leq(k+1)\ell-2$,
\begin{eqnarray*}
\frac{z(n+2)-z(n+1)}{z(n+1)-z(n)}=L^{1/\ell}.
\end{eqnarray*}
Also
\begin{eqnarray*}
z(k\ell)-z(k\ell-1)=L^{(\ell-1)/\ell}\frac{L^{1/\ell}-1}{L-1}(v(k)-v(k-1))\end{eqnarray*}
and
\begin{eqnarray*}
z(k\ell+1)-z(k\ell)=\frac{L^{1/\ell}-1}{L-1}(v(k+1)-v(k))\leq\frac{L^{1/\ell}-1}{L-1}L(v(k)-v(k-1)),
\end{eqnarray*}
so that the ratio
\begin{eqnarray*}
\frac{z(k\ell+1)-z(k\ell)}{z(k\ell)-z(k\ell-1)}\leq L^{1/\ell}
\end{eqnarray*}
as well. This shows that $z(n+2)-z(n+1)\leq L^{1/\ell}(z(n+1)-z(n))$ for all $n$. Since $L^{1/\ell}<2$ and $z(n+1)-z(n)$ is bounded below, there exists a large constant $C$ such that, when $z(n)$ is replaced by $w(n)=\lfloor Cz(n)\rfloor$, the inequality $w(n+2)-w(n+1)\leq 2(w(n+1)-w(n))$ remains valid. This also makes $w(n+1)-w(n)\geq 2$. Since $v(n)=O(L^n)$, $w(n)=O(L^{n/\ell})$ and $L^{1/\ell}<2$. Clearly, $w$ has the same growth type as $v$. And if $v(n+1)-v(n)$ tends to $+\infty$, so does $w(n+1)-w(n)$. Substracting a constant makes $w(0)=1$.
\end{pf}

\subsection{Further requirements on pieces}
\label{req}

\begin{nota}
\label{notavolumes}
For a piece $P$, let $t_{P}$ (resp. $T_{P}$) denote the minimum
(resp. maximum) of the function distance to $\partial^{-}P$ restricted to
$\partial^{+}P$. For $k\leq T_{P}$, let $U_{P,k}$ denote the $k$-tubular neighborhood of $\partial^{-}P$ and
\begin{eqnarray*}
v_{P}(k)=vol(U_{P,k}),\quad v'_{P}(k)=v_{P}(k)-v_{P}(k-1).
\end{eqnarray*}
\end{nota}

\begin{prop}
\label{pieces}
Let $Q_j$ be a sequence of compact manifolds with boundary. Assume that
\begin{itemize}
  \item $\partial Q_j$ is split into two open and closed subsets
    $\partial^{-}Q_j$ and $\partial^{+}Q_j$;
  \item $\partial^{-}Q_{j+1}$ is diffeomorphic to $\partial^{+}Q_j$.
\end{itemize}
Then there exist integers $\ell$, $h$, $H$, sequences of integers $t_j$,
$u_j$, $U_j$, $d_{j}$ and Riemannian metrics on pieces $Q_j$, $R_j$, $K$, $HS$, $J$ such that
\begin{enumerate}
  \item For all pieces $P$, the maximal distance of a point of $P$ to
    $\partial^{-}P$ is achieved on $\partial^{+}P$. In other words, it is
    equal to $T_{P}$.
  \item $\frac{1}{3}\ell t_{j}\leq t_{Q_{j}}\leq
    T_{Q_{j}}\leq \ell t_{j}$.
  \item For all other pieces $P$, $\frac{1}{3}\ell\leq t_{P}\leq
    T_{P}\leq\ell$.
  \item $diameter(\partial^{-}Q_{j})\leq d_{j}$.
  \item All pieces $P$ carry a marked point $y_{P}\in\partial^{-}P$. When
    a piece $P'$ is glued on top of $P$, $d(y_{P},y_{P'})\leq \ell$
    (resp. $\ell t_{j}$ if $P=Q_{j}$), unless
    $P=R_{j}$ and $P'$ is of type $K$, $HS$ of $J$. In that case, $d(y_{P},y_{P'})\leq d_{j}$.
  \item For all pieces $P=K$, $HS$, $J$, $h\leq\min v'_{P}\leq \max v'_{P}\leq H$.
  \item $\max v'_{Q_j}\leq U_j$.
  \item $\max v'_{R_j}\leq u_{j}\leq U_{j}$.
  \item If $\partial^{+}Q_{j}$ and $\partial^{-}Q_{j}$ are diffeomorphic,
    then they are isometric, by an isometry that maps $y_{j}$ to $y_{j+1}$,
    and $u_{j+1}= u_{j}$.
  \item All pieces have bounded geometry and product metric near the boundary.
\end{enumerate}
$t_j$,
$u_j$, $d_{j}$ are respectively called the \emph{height, volume} and \emph{diameter parameters}.
\end{prop}

\subsection{Proof of Proposition \ref{pieces}}
\label{proofpieces}

The cases of cylinder $K$ and half-sphere $HS$ are easy. For $Q_{j}$, we shall start with some initial metric satisfying weak requirements, and modify it in two steps,
\begin{enumerate}
  \item glue product manifolds $[-T,0]\times\partial P$ along the
    boundary, equipped with warped product metrics modelled on
    hyperbolic cusps.
  \item rescale so that the metric has bounded geometry and the boundary gets back to its original size.
\end{enumerate}
For $R_{j}$ and the join $J$, rescaling, and thus warping, is unneeded:
thickening the boundary with direct product metrics is sufficient to achieve (1) and (3).

\subsubsection{Initial metric on $Q_{j}$}

Choose a point $y_{Q_{j}}\in\partial^{-}Q_j$. The only constraint is the
following : if two consecutive manifolds $\partial^{-}Q_j$ and
$\partial^{-}Q_{j+1}$ are diffeomorphic, pick one such diffeomorphism
$\phi_{j}:\partial^{-}Q_j \to\partial^{-}Q_{j+1}$ and assume that $y_{Q_{j+1}}=\phi_{j}(y_{Q_{j}})$.

Pick a Riemannian metric of bounded geometry on each of the manifolds
$\partial^{-}Q_j$. The only constraint is the following : if two
consecutive manifolds $\partial^{-}Q_j$ and $\partial^{-}Q_{j+1}$ are
diffeomorphic, pick isometric metrics (i.e. mapped to each other by the
chosen diffeomorphism $\phi_{j}$). Modify it slightly so that it is flat on some ball of radius $3$. Extend the resulting metric on
$$\partial^{-}Q_{j}\coprod\partial^{-}Q_{j+1}=\partial^{-}Q_{j}\cup \partial^{+}Q_j =\partial Q_j $$
to a product metric on some collar neighborhood of $\partial Q_j$ in
$Q_j$. Extend it arbitrarily to a Riemannian metric $m_{j}$ on $Q_j$. Let
\begin{eqnarray*}
\lambda_j =\max\{(\mathrm{injectivity~radius})^{-1},\sqrt{\mathrm{Max~sectional~curvature}}\}
\end{eqnarray*}
be the scaling factor needed to turn $m_{j}$ into a metric of bounded geometry.

\subsubsection{Initial metric on the join $J$}

Start with the Euclidean metric on $\R^{n}$. By a conformal change supported in the $3$-ball centered at the origin, turn it into a complete metric of revolution on $\R^{n}\setminus B(0,1)$ which, in the $2$-ball, is isometric to a direct product $[0,1]\times S^{n-1}$. Call this a handle. Now start with the product metric on $[-10,10]\times S^{n-1}$. Without changing the boundary, modify it to make it flat in a ball of radius $3$, then surge in the handle.

\subsubsection{Initial metric on $R_{j}$}

Use the previously chosen metric on $\partial^{-}R_{j}=\partial^{+}
Q_j$. The product $[-10,10]\times \partial^{-}R_{j}$ has bounded
geometry and contains a flat ball of radius $3$. Surge in a handle in the flat part, to produce a new boundary component, isometric to a unit round sphere.

\subsubsection{Thickening the boundary}

For $P=J$ or $R_j$, thickening merely means gluing in a Riemannian product $[-T,0]\times\partial P$.

For $Q_j$, let us proceed in two steps. First glue to $Q_j$ a Riemannian product $[-T,0]\times\partial Q_j$, leading to a metric $m_{j,T}$ on $Q_j$ in which the $T$-tubular neighborhood of the boundary is a product.

Fix once and for all a smooth positive nondecreasing function $f:\R\to\R$ such that
\begin{itemize}
  \item $f(t)=1$ for $t>0$;
  \item $f(t)=e^t$ for $-2<t<-1$;
  \item $f(t)=e^{-3}$ for $t<-4$.
\end{itemize}
One can assume that $f$ is convex and satisfies $f''\geq f$ on $(-\infty,-1]$. Define $f_T : [-T,0]\to\R$ by
\begin{itemize}
  \item $f_{T}(t)=f$ for $t>-2$;
  \item $f_{T}(t)=e^t$ for $-T+2<t<-1$;
  \item $f_{T}(t)=e^{-T+5}f(t+T-5)$ for $-T\leq t<-T+2$.
\end{itemize}
Change the metric $m_{j,T}$ in the $T$-tubular neighborhood of the boundary from a product metric $dt^2 +\tilde{g}$ to a warped product $g_{j,T}=dt^2 +f_{T}(t)^2 \tilde{g}$. This new metric is still a product in the $1$-neighborhood of the boundary.

\subsubsection{Controlling the heights of $J$ and $R_j$}

We claim that if $T$ is large enough, properties (1) and (3) are satisfied.

The argument applies indifferently to $P=J$ or $R_j$ and in the latter case,
$T$ does not depend on $j$. Let $D$ denote the diameter of the region in
the handle where the metric is not a product. If $T>D+20$, a point $q$ of
$P$ at maximal distance from $\partial^{-}P$ must belong to the
$T$-neighborhood $W^{+}$ of $\partial^{+}P$. Let $p$ be a point of
$\partial^{-}P$ closest to $q$. Let $\gamma$ be a minimizing geodesic from
$p$ to $q$. In $W^{+}$, the derivative $\gamma'$ makes a constant angle with the
$\partial^{+}P$ factors, thus $\gamma'$ points towards $\partial^{+}P$ at
$q$. If $q\notin\partial^{+}P$, one could move $q$ towards $\partial^{+}P$
and increase distance to $\partial^{-}P$, contradiction. We conclude that
$q\in\partial^{+}P$, this is (1) for $J$ and $R_j$.

Every point of $P$ sits at a distance from $\partial^{-}P$ at most
$2T+20+D$. Every point of $\partial^{+}P$ sits at a distance from
$\partial^{-}P$ at least $T$. So since $T>D+20$, $T\leq t_{P}\leq T_{P}\leq
3T$. With $\ell=3T$, this is (3). Note that one still may enlarge $T$ (and
thus $\ell$), provided it does not depend on $j$.

\subsubsection{Heights in warped products}

We shall need the following Lemma.

\begin{lem}
\label{hyper}
Let $(\tilde{M},\tilde{g})$ be a complete Riemannian manifold. Let $M=[-T,0]\times\tilde{M}$ be equipped with the warped product metric $dt^2 +f(t)^2 \tilde{g}$. Let $t_{0}\in[-T,0]$. Let $m$, $m'\in\tilde{M}$. Let $s\mapsto\gamma(s)=(t(s),\tilde{\gamma}(s))$ denote a minimizing geodesic from $(t_{0},m)$ to $(0,m')$ in $M$. Then
\begin{enumerate}
  \item $\tilde{\gamma}$ is a minimizing geodesic from $m$ to $m'$ in $\tilde{M}$.
  \item If $d_{\tilde{g}}(m,m')<e^{-1-t_{0}}-3$, then $s\mapsto t(s)$ is monotone with derivative $t'(0)>0$.
\end{enumerate}
\end{lem}

\begin{pf}
The second fundamental form of the hypersurface $\{t\}\times\tilde{M}$ is $f'(t)\tilde{g}$. Thus if $u=a\frac{\partial}{\partial t}+v$ is a vectorfield on $M$, the tangential component of the Levi-Civita connection is
\begin{eqnarray*}
\nabla^{tan}_{u}u=\tilde{\nabla}_{v}v-2a\frac{f'}{f}v.
\end{eqnarray*}
If $\gamma$ is a constant speed geodesic, $\nabla^{tan}_{\gamma'}\gamma' =0$, $\tilde{\nabla}^{tan}_{\tilde{\gamma}'}\tilde{\gamma}' $ is colinear to $\tilde{\gamma}'$, hence $\tilde{\gamma}$ is a reparametrization of a geodesic. Note that the speed of $\tilde{\gamma}$ does not vanish, unless it vanishes identically. To compute the length of $\gamma$, one may restrict to the immersed submanifold $[-T,0]\times\tilde{\gamma}$, i.e. assume that $\tilde{M}=[0,L]$, $L=length(\tilde{\gamma})$. If $t_0$ and $T$ are fixed, $length(\gamma)$ is a function of $L$ only. This function is increasing. Indeed, if $L'>L$, there exists $s_{0}$ such that $\tilde{\gamma}_{L'}(s_{0})=L$. Replacing the arc of $\gamma_{L'}$ from $\gamma_{L'}(s_{0})$ to $(0,L')$ with a segment of the form $s\mapsto (s,L)$ to obtain a curve from $(t_0 ,0)$ to $(0,L)$ reduces length, showing that $length(\gamma_{L'})>length(\gamma_{L})$. Therefore, if $\gamma$ is length minimizing, so is $\tilde{\gamma}$.

Again, in order to study the sign of $t'(0)$, one may assume that $\tilde{M}=[0,L]$. We first reason on $[-T,-1]\times[0,L]$. There, curvature is nonpositive, this draws geodesics backwards, and makes certain shortest geodesics have $t'(0)<0$. Fortunately, curvature stays $\geq -1$. To get estimates, it suffices to treat the case when curvature equals $-1$ everywhere, i.e. to study geodesics in the hyperbolic plane. The change of coordinates $(t,x)\mapsto (x,e^{t})$ maps to the upper half plane model. Thus hypersurfaces $\{t\}\times[0,L]$ are pieces of horocycles. If a geodesic starts tangentially to a horocycle and reaches a point on the parallel horocycle at distance $-1-t_0$, whose abscissa is $x$, then $x^2 +1=e^{-2(1+t_{0})}$ (see figure). Thus if $x<\sqrt{e^{-2(1+t_{0})}-1}$, the geodesic from $(t_0 ,0)$ to $(-1,x)$ has $t'(0)>0$. If follows that if $L<\sqrt{e^{-2(1+t_{0})}-1}-2$, the geodesic from $(t_0 ,0)$ to $(0,L)$ has $t'(0)>0$.
\end{pf}

\subsubsection{Controlling the height of $Q_j$}

We claim that if $T=T(j)$ is large enough, properties (1) and (2) are satisfied.

Let $D(j)=diameter(Q_j ,m_j )$. The diameter of $(Q_j ,g_{j,T})$ lies between
$2T$ and $2T+3D$. If $T>2D$, a point $q$ at maximal distance from
$\partial^{-}Q_{j}$ must sit in $W^{+}$.

Let $q$ have coordinates $(t_0 ,m)$ in $W^{+}$. Then $t_0
+T<D$. According to Lemma \ref{hyper}, since $t_0 <-\log(D(j))$ (roughly),
all minimizing geodesics from $p$ to $q$ point away from the
boundary. Pulling $q$ forwards, i.e. towards $\partial^{+}Q_j$ should
allow to increase distance from $\partial^{-}Q_{j}$. We conclude that
$q\in\partial^{+}Q_j$. This proves property (1).

Every point of $Q_{j}$ sits at a distance from $\partial^{-}P$ at most
$2T+D(j)$. Every point of $\partial^{+}P$ sits at a distance from
$\partial^{-}P$ at least $T$. So as soon as $T>2D(j)$,
\begin{eqnarray*}
T\leq t_{Q_{j}}\leq T_{Q_{j}}+2D\leq T+4D(j)\leq 3T.
\end{eqnarray*}

Now rescale the metric on $Q_{j}$, i.e., replace metric $g_{j,T}$ by
$g'_{j,T}=e^{2T-4}g_{j,T}$. The metric induced by $g'_{j,T}$ on the
boundary does not depend on $T$, it has bounded geometry. The exponential
warping does not spoil curvature bounds. Furthermore,
since $\lambda_{j}^{-2}m_{j}$ is a bounded geometry metric on
$Q_j$, $g_{j,T}$ has bounded geometry as soon as $e^{T-2}\geq\lambda_j$. So
property (10) holds for the rescaled $Q_{j}$ if we take
$T=T(j)=\max\{2D(j),\log(\lambda_{j})\}$. Also, the scale invariant inequality
\begin{eqnarray*}
t_{Q_{j}}\leq T_{Q_{j}}+2diameter(\partial Q_{j})\leq 3t_{Q_{j}}
\end{eqnarray*}
still holds. Finally, one may choose an integer $t_{j}$ such that $\frac{1}{3}\ell t_{j}\leq
t_{Q_{j}}\leq T_{Q_{j}}+2diameter(\partial Q_{j})\leq \ell t_{j}$, this is (2).

\subsubsection{Fixing the diameter parameter}

When small pieces $P'$ of type $J$, $K$, $HS$ are glued on top of each
other, $d(y_{P},y_{P'})$ is bounded, so one can assume that
$d(y_{P},y_{P'})\leq \ell$. When a piece $P'$ of type $Q$ or $R$ is glued on top of a piece of type
$R$, $y_{P'}$ is on top of $y_{P}$, thus $d(y_{P},y_{P'})\leq \ell$.
When $R_{j}$ is glued on top of $Q_{j}$, $d(y_{Q_{j}},y_{Q_{j}})\leq
T_{Q_{j}}+2diameter(\partial Q_{j})\leq\ell t_{j}$. The only bad case
happens when a piece $P'$ of type $J$, $K$, $HS$ is glued on top of
$R_{j}$. In this case, we simply define
$d_{j}=\max\{diameter(\partial^{-}(Q_{j}),d(y_{R_{j}},y_{P'}\}$, so
properties (4) and (5) hold.

\subsubsection{Fixing volume parameters}

For small pieces, there are finitely many slice volumes $v'_{P}(k)$, which
are bounded from below by some $h$ and above by some $H$, and (6) holds.
Define $u_{j}=\max v'_{R_{j}}$ and $U_{j}=\max\{u_{j},\max v'_{Q_{j}}\}$,
so that (7) and (8) hold.

Properties (9à and (10) hold by construction. This completes the proof of
Proposition \ref{pieces}.

\subsection{The discrete growth function}

\begin{defi}
\label{defz}
Let $T$ be an admissible rooted tree. Let $v:\N\to\N$ denote its growth. Glue together pieces according to the pattern given by $T$ and get a Riemannian manifold $R_{T}$. Define a function $r:R_{T}\to\N$ as follows. If $P=Q_j$ and $x\in Q_j$, let $r(x)=\lfloor d(x,\partial^{-}Q_j)\rfloor+n_j \ell$. If $P$ is any other type of piece, attached at a vertex of $T$ of level $n$, and $x\in P$, let $r(x)=\lfloor d(x,\partial^{-}P)\rfloor+n\ell$. We define the associated {\em discrete growth function} $z$ as follows: for $n\in\N$, $z(n)=vol(\{x\in R_{T}\,|\,r(x)\leq n\}$.
\end{defi}

\begin{lem}
\label{z}
Here are bounds for the discrete growth function $z$ associated to a given function $v$ satisfying the assumptions of Lemma \ref{tree}. If $\ell n_j \leq n <\ell n_{j+1}$,
$$(v(n)-v(n-1)-1)h\leq z(n)-z(n-1)\leq H(v(n)-v(n-1))+U_{j}.$$
\end{lem}

\begin{pf}
The set $\{x\in R_{T}\,|\,r(x)=n\}$ is a union of $v(n)-v(n-1)$ slices taken from various types of pieces. The first of these pieces (in lexicographical order), is either a $Q_j$ or a $R_j$, therefore the volume of the slice is at most $U_{j}$ or $u_{j}$, and both are less than $U_{j}$. In all other pieces, the volumes of slices are at least $h$ and at most $H$. Adding up theses volumes yields the volume $z(n)-z(n-1)$ of $\{x\in R_{T}\,|\,r(x)=n\}$.
\end{pf}

\begin{prop}
\label{prall}
Let $v:\N\to\N$ be a function that satisfies the assumptions of Lemma \ref{tree}, i.e.
\begin{itemize}
\item $v(0)=1$.
\item for all $n\in\N$, $2\leq v(n+2)-v(n+1)\leq 2(v(n+1)-v(n))$.
\item $v(n)=O(\lambda^n)$ for some $\lambda<2$.
\end{itemize}
Let $t_j$, $u_j$, $d_{j}$ be the parameters of the pieces $Q_j$, as provided by Lemma \ref{pieces}. Assume that
\begin{itemize}
  \item either $\lim_{n\to \infty}v(n+1)-v(n)=+\infty$;
  \item or $u_j$ is constant.
\end{itemize}
Then there exists an increasing sequence $n_j$ such that
\begin{enumerate}
  \item $n_{j}\geq d_{j}$.
  \item the subset $\displaystyle S=\bigcup_{j}[n_j ,n_j +t_j -1]$ has vanishing lower density;
  \item the discrete growth function $z$ of the corresponding Riemannian manifold $\displaystyle R_{T_{S,v}}$ has the same growth type as $v$.
\end{enumerate}
\end{prop}

\begin{pf}
Assume first that $v(n+1)-v(n)$ tends to infinity. $n_j$ will be chosen inductively. Let us collect specifications for $n_j$. By assumption, there exists $r_j$ such that $v'(n):=v(n)-v(n-1)\geq \max\{h,U_{j}\}$ for $n\geq r_j$. So we require $n_j \geq r_j$. According to Lemma \ref{z}, for $n_j \leq n <n_{j+1}$,
\begin{eqnarray*}
h(v'(n)-1) \leq z'(n)\leq Hv'(n)+U_{j}.
\end{eqnarray*}
This implies that
\begin{eqnarray*}
(h-1)v'(n)\leq z'(n)\leq (H+1)v'(n).
\end{eqnarray*}
The other specification is that the union $S$ of intervals $[n_j ,n_j
+t_j]$ have vanishing lower density. This is achieved by requiring that
$n_{j}\geq j(n_{j-1}+t_{j-1})$. Thus we take $n_j =\max\{d_{j},t_j ,r_j
,j(n_{j-1}+t_{j-1})\}$.

Assume next that $u_j =u$ does not depend on $j$. We
first construct an admissible tree $T_{0}$ whose branches have depth $1$,
except for the trunk. In other words, $T_{0}$ consists in a ray, the trunk,
with one edge glued at a trunk vertex $x_{k}$ unless $n_{j}\leq
k<n_{j}+t_{j}$. Let $z_{0}$ denote the discrete growth function of the
corresponding manifold $R_{0}$. We pick $n_{j}$ inductively in such a way
that $z_{0}(n)\leq 4u^{2}n$. Assume that $n_{j-1}$ has been defined, and
that $z_{0}(n)\leq 4u^{2}n$ for all $n\leq n_{j-1}+t_{j-1}$. Lemma
\ref{stretch} below, applied with $A=\max v'_{Q_{j}}$, $B=2u$, $C=u$,
$a=n_{j-1}+t_{j-1}$, $b=t_{j}$, $v(n)=n$, provides us with $R$ such that if
we take $n_{j}=n_{j-1}+t_{j-1}+R$, then $z_{0}(n)\leq B^{2}n=4u^{2}n$ for
$n\leq n_{j}+t_{j}$ (this construction is taken from \cite{G1}). Since $R$
can be chosen arbitrarily large, there is no obstacle to let the set $S$
have vanishing lower density and to achieve $n_{j}\geq d_{j}$. Lemma \ref{tree} upgrades $T_0$ into an admissible tree $T$ with volume growth $v$. The corresponding Riemannian manifold $R$ has discrete growth function $z$. Let $v_0$ denote the volume growth of $T_0$. Note that $v'_{0}(n)=v_0 (n+1)-v_0 (n)$ takes only two values, $1$ or $2$. Lemma \ref{z} implies that if $\ell n_j \leq n <\ell n_{j+1}$,
$$h(v'(n)-v'_{0}(n))\leq z'(n)-z'_{0}(n)\leq Hv'(n).$$
Integrating yields
\begin{eqnarray*}
z_{0}(n)+h(v(n)-v_{0}(n))\leq z(n)\leq z_{0}(n)+Hv(n).
\end{eqnarray*}
Since $v$ grows at least linearly, and $z_0$ and $v_0$ at most linearly, this shows that $z$ has the same growth type as $v$.
\end{pf}

\begin{lem}
\label{stretch} Let $a>0$, $b>0$, $A>B>C>1$. Let $v$ be a
nondecreasing function on $\R_+$. There exist arbitrarily large
$R=R(a,b,A,B,C,v)$ such that if a nondecreasing function $z$ on
$[0,a+R+b]$ satisfies
\begin{enumerate}
  \item $B^{-1}v(B^{-1}x)\leq z(x)\leq Bv(Bx)$ on $[0,a]$,
  \item $C^{-1}v'\leq z'\leq Cv'$ on $[a,a+R]$,
  \item $A^{-1}v'\leq z'\leq Av'$ on $[a+R,a+R+b]$,
\end{enumerate}
then
\begin{eqnarray*}
\forall x\in [0,a+R+b],\quad B^{-1}v(B^{-1}x)\leq z(x)\leq Bv(Bx).
\end{eqnarray*}
\end{lem}

\begin{pf}
Clearly, since $C\leq B$, if $x\leq a+R$, $B^{-1}v(B^{-1}x)\leq
z(x)\leq Bv(Bx)$. Let $x>a+R$. Then
\begin{eqnarray*}
z(x)&=&z(a)+\int_{a}^{a+R}z'(t)\,dt+\int_{a+R}^{x}z'(t)\,dt\\
&\leq& Bv(Ba)+C\int_{a}^{a+R}v'(t)\,dt+A \int_{a+R}^{x}v'(t)\,dt\\
&\leq&Bv(Ba)+Cv(a+R)+A(v(a+R+b)-v(a+R))\\
&=:& f(R).
\end{eqnarray*}
Also,
\begin{eqnarray*}
z(x)
&\geq& B^{-1}v(B^{-1}a)+C^{-1}\int_{a}^{a+R}v'(t)\,dt+A^{-1}\int_{a+R}^{x}v'(t)\,dt\\
&\geq&B^{-1}v(B^{-1}a)-C^{-1}v(a)+C^{-1}v(a+R)\\
&=:& g(R).
\end{eqnarray*}
If $\displaystyle \liminf_{R\to\infty}\frac{v(a+R+b)}{v(a+R)}=1$,
then there exists a sequence $R_j$ tending to $\infty$ such that
\begin{eqnarray*}
\lim_{j\to \infty}\frac{f(R_j)}{v(a+R_j)}=C,\quad \lim_{j\to
\infty}\frac{g(R_j)}{v(a+R_j +b)}=C^{-1}.
\end{eqnarray*}
Thus for $j$ large enough, $\displaystyle
\frac{f(R_j)}{v(a+R_j)}\leq B$ and $\displaystyle
\frac{g(R_j)}{v(a+R_j +b)}\geq B^{-1}$. Pick such a $j$, then, for $x\geq a+R_j$,
\begin{eqnarray*}
z(x)\leq f(R_j)\leq Bv(a+R_j)\leq Bv(x),\quad z(x)\geq g(R_j)\geq
B^{-1}v(a+R_j +b)\geq B^{-1}v(x).
\end{eqnarray*}
Otherwise, $\displaystyle
\liminf_{R\to\infty}\frac{v(a+R+b)}{v(a+R)}>1$.  Then there exists
$\lambda>1$ and $r_0$ such that $r\geq r_0 \Rightarrow
v(r+b)\geq\lambda v(r)$. For $r\geq r_0$,
\begin{eqnarray*}
v(Br)\geq\lambda^{(B-C)r/b}v(Cr),
\end{eqnarray*}
i.e. $v(Br)/v(Cr)$ tends to $+\infty$. In particular, $v$ tends to
$+\infty$. Take $R$ large enough so that $a+R+b\leq C(a+R)$. Then
\begin{eqnarray*}
\frac{z(x)}{Bv(Bx)}&\leq&\frac{f(R)}{Bv(B(a+R))}\\
&\leq& \frac{v(Ba)}{v(B(a+R))}+\frac{A(v(a+R+b)}{Bv(B(a+R))}\\
&\leq& o(1)+\frac{A}{B}\frac{v(C(a+R))}{v(B(a+R))}
\end{eqnarray*}
which tends to $0$, so is $\leq 1$ for $R$ large enough.
Similarly,
\begin{eqnarray*}
\frac{z(x)}{B^{-1}v(B^{-1}x)}&\geq&\frac{g(R)}{B^{-1}v(B^{-1}(a+R+b))}\\
&\geq& \frac{B^{-1}v(B^{-1}a)-C^{-1}v(a)}{B^{-1}v(B^{-1}(a+R+b))}+\frac{C^{-1}v(a+R)}{B^{-1}v(B^{-1}(a+R+b))}\\
&\geq& o(1)+\frac{B}{C}
\end{eqnarray*}
which is $\geq 1$ for $R$ large enough.
\end{pf}

\medskip

Note that the assumption $\lim_{n\to\infty}v(n+1)-v(n)=+\infty$ of Proposition \ref{prall} is not that restrictive. Up to changing of representative of a growth type, it follows from the weaker assumption $\lim_{n\to\infty}\frac{v(n)}{n}=+\infty$ of Theorem \ref{infinite}.

\begin{lem}
\label{suplinear}
Let $v:\N\to\N$ be a non decreasing function satisfying for all $n\in \N$,
\begin{eqnarray*}
v(n+2)-v(n+1)\leq L(v(n+1)-v(n)).
\end{eqnarray*}
Assume that
\begin{eqnarray*}
\lim_{n\to\infty}\frac{v(n)}{n}=+\infty.
\end{eqnarray*}
Then there exists a function $w:\N\to\R$, having the same growth type as $v$, such that
\begin{enumerate}
  \item for all $n\in \N$,
\begin{eqnarray*}
\frac{1}{L}\leq w(n+2)-w(n+1)\leq L(w(n+1)-w(n)).
\end{eqnarray*}
\item $\lim_{n\to\infty}w(n+1)-w(n)=+\infty$.
\end{enumerate}
\end{lem}

\begin{pf}
Let
\begin{eqnarray*}
Y=\{(x,y)\in\R^2 \,;\,x\in\N,\,y\geq v(x)\}
\end{eqnarray*}
denote the epigraph of $v$, let $Z\subset\R^2$ be its convex hull, and
\begin{eqnarray*}
u(x)=\min\{y\in\R\,;\,(x,y)\in Z\}.
\end{eqnarray*}
By construction, $u\leq v$, and $u$ is convex. By assumption, for every
line $L$ through the origin with positive and finite slope, the part of $Y$
below $L$ is compact. Therefore the part of $Z$ below $L$ is compact as well. This shows that
\begin{eqnarray*}
\lim_{x\to\infty}u'(x)=+\infty,
\end{eqnarray*}
and thus
\begin{eqnarray*}
\lim_{n\to\infty}u(n+1)-u(n)=+\infty.
\end{eqnarray*}
By construction, $u$ is piecewise linear and its derivative changes only at integers $n$ such that $u(n)=v(n)$. At such points,
\begin{eqnarray*}
u(n+1)-u(n)\leq v(n+1)-v(n), \quad u(n)-u(n-1)\geq v(n)-v(n-1),
\end{eqnarray*}
thus $u(n+1)-u(n)\leq L(u(n)-u(n-1))$. At other points, $u(n+1)-u(n)=u(n)-u(n-1)\leq L(u(n)-u(n-1))$.

Set $w(n)=u(n)+v(n)$. Then for all $n$, $w(n+1)-w(n)\leq L(w(n)-w(n-1))$,
$v(n)\leq w(n)\leq 2w(n)$. Furthermore,
$\lim_{n\to\infty}w(n+1)-w(n)=+\infty$. Adding a constant and further
changing $w$ at finitely many places allows to have $w(n+1)-w(n)\geq 2$ for
all $n$.
\end{pf}

\subsection{End of the proof of Theorem \ref{infinite}}
\label{proofinfinite}

Let $v$ be a given bgd-function.

If $M$ has finite topological type, use an exhaustion of $M$ into $M_j$'s with boundaries diffeomorphic to a fixed compact manifold $V$. Proposition \ref{pieces} provides us with pieces, and in particular, Riemannian metrics on $Q_j =M_{j+1}\setminus M_{j}$, with isometric boundaries and constant  volume parameters $u_j$. Proposition \ref{prall} shows how to adjust remaining parameters (a sequence $(n_j)$) so that the discrete growth function $z$ of the resulting Riemannian manifold $R$ is equivalent to $v$.

If $M$ has infinite topological type, we assume that $\lim v(n)/n=+\infty$. Lemma \ref{suplinear} even allows us to assume that $\lim v(n+1)-v(n)=+\infty$, so that Proposition \ref{prall} applies again.

In both cases, there remains to relate $z$ to the true growth function of
$R$, i.e. $w(n)=vol(B(o,n))$. We first relate the integer valued function
$r$ to the Riemannian distance to $o$. Let $x\in R$ belong to some piece
$P$ at level $n$, so that $r(x)=\lfloor d(x,\partial^{-}P)\rfloor+n\ell$
(resp. $+n_{i} \ell$ if $P=Q_{i}$ and $n_{i}\leq n<n_{i}+t_{i}$).

Connect successive marked points $y_{0}=o$,
$y_{1},\ldots,y_{k}\in\partial^{-}P$. Let $y_{k+1}$ be the point of
$\partial^{-}P$ which is closest to $x$. Connect $y_{k}$ to $y_{k+1}$ in
$\partial^{-}P$ and $y_{k+1}$ to $x$ in $P$ by minimizing geodesics. The constructed path yields the estimate
\begin{eqnarray*}
d(o,x)\leq\sum_{i=0}^{k}d(y_{i},y_{i+1})+d(y_{k+1},x).
\end{eqnarray*}
According to Proposition \ref{pieces},
$d(y_{i},y_{i+1})\leq\ell$ unless $y_{i}$ belongs to a piece $P_{i}$ of type
$R$ and $y_{i+1}$ to a piece of type $K$, $HS$ or $J$. This may happen
for at most one value of $i$ in $\{1,\ldots,k\}$, and in that case,
$d(y_{i},y_{i+1})\leq diameter(\partial^{-}P_{i})+\ell$. $P_{i}$ belongs to
a pile of $R$'s glued to a $Q_{j}$, and
$diameter(\partial^{-}P_{i})=diameter(\partial^{-}Q_{j})\leq n_{j}\ell$.

If $P$ is of type $Q$, e.g. $P=Q_{m}$, then $j\leq m$,
\begin{eqnarray*}
\sum_{i=0}^{k-1}d(y_{i},y_{i+1})\leq 2n_{m}\ell.
\end{eqnarray*}
Furthermore,
$d(y_{k},y_{k+1})\leq diameter(\partial^{-}Q_{m})\leq n_{m}\ell$, thus
\begin{eqnarray*}
\sum_{i=0}^{k}d(y_{i},y_{i+1})\leq 3n_{m}\ell,
\end{eqnarray*}
and
\begin{eqnarray*}
d(o,x)\leq 3n_{m}\ell+(r(x)-n_{m}\ell)\leq 3r(x).
\end{eqnarray*}

Otherwise,
\begin{eqnarray*}
\sum_{i=0}^{k-1}d(y_{i},y_{i+1})\leq 2n\ell.
\end{eqnarray*}
Furthermore, $d(y_{k},y_{k+1})\leq\ell$, and
\begin{eqnarray*}
d(o,x)\leq (2n+1)\ell+(r(x)-n\ell)\leq \frac{n+1}{n}r(x)\leq 3r(x).
\end{eqnarray*}

Conversely, let $\gamma$ be a minimal geodesic segment from $o$ to $x$. It
passes through $n$ pieces (where $Q_{j}$'s are counted with multiplicity
$t_{j}$). Let $s_0 =o$ and let $s_1,s_2 ,\ldots,s_k$ denote the values of
$s$ such that $\gamma(s)$ belongs to the lower boundary $\partial^{-}P$ of
some piece. By construction, $d(\gamma(s_i),\gamma(s_{i+1}))\geq\ell/3$
(resp. $\geq\ell t_j /3$, depending on the type of piece). Also
$d(\gamma(s_k),x)\geq d(x,\partial^{-}P)$. Summing up yields
\begin{eqnarray*}
d(o,x)\geq \frac{1}{3}n_{m}\ell+d(x,\partial^{-}P)\geq
\frac{1}{3}n_{m}\ell+(r(x)-n_{m}\ell)\geq \frac{1}{2}r(x),
\end{eqnarray*}
if $P=Q_{m}$, and
\begin{eqnarray*}
d(o,x)\geq \frac{1}{3}n\ell+d(x,\partial^{-}P)\geq
\frac{1}{3}n\ell+(r(x)-n\ell)\geq \frac{1}{3}r(x),
\end{eqnarray*}
otherwise. This shows that
\begin{eqnarray*}
\{x\in R\,|\,r(x)\leq \frac{n}{3}\}\subset B(o,n)\subset \{x\in R\,|\,r(x)\leq 3n\},
\end{eqnarray*}
and $z(\frac{n}{3})\leq vol(B(o,n)) \leq z(3n)$. One concludes that the constructed bounded geometry manifold has volume growth equivalent to $v$.

\bibliographystyle{plain}

\end{document}